\documentclass{article}%
\usepackage{amssymb}
\usepackage{amsfonts}
\usepackage{amsmath}
\usepackage{graphicx}%
\begin{document}

\title{A Novel Solution to the Frenet-Serret Equations}
\author{Anthony A. Ruffa\\Naval Undersea Warfare Center Division\\1176 Howell Street\\Newport, RI 02841}
\maketitle

\begin{abstract}
A set of equations is developed to describe a curve in space given the
curvature $\kappa$ and the angle of rotation $\theta$ of the osculating plane.
\ The set of equations has a solution (in terms of $\kappa$ and $\theta$) that
indirectly solves the Frenet-Serret equations, with a unique value of $\theta$
for each specified value of $\tau$. \ Explicit solutions can be generated for
constant $\theta$. \ The equations break down when the tangent vector aligns
to one of the unit coordinate vectors, requiring a reorientation of the local
coordinate system.

\end{abstract}

\section{Introduction}

Given the curvature $\kappa$ and torsion $\tau$, the Frenet-Serret
equations$^{1}$ describe a curve in space parameterized by the arc length $s$:%

\begin{align}
\frac{d\boldsymbol{T}}{ds}  &  =\kappa\boldsymbol{N};\\
\frac{d\boldsymbol{N}}{ds}  &  =-\kappa\boldsymbol{T+}\tau\boldsymbol{B};\\
\frac{d\boldsymbol{B}}{ds}  &  =-\tau\boldsymbol{N};\\
\frac{d\boldsymbol{R}}{ds}  &  =\boldsymbol{T}.
\end{align}

Here $\boldsymbol{R}$, $\boldsymbol{T}$, $\boldsymbol{N}$, and $\boldsymbol{B}%
$ are the position, tangent, normal, and binormal vectors, respectively.
\ These equations have no explicit solution (in terms of $\kappa$ and $\tau$)
for the general case, although solutions for special cases exist$^{2}$.

It is shown here that a set of equations can be developed to describe a curve
in space given the curvature $\kappa$ and the angle of rotation $\theta$ of
the osculating plane. \ The set of equations has a solution (in terms of
$\kappa$ and $\theta$) that indirectly solves the Frenet-Serret equations, and
has a unique $\theta$ for every value of $\tau$.

Many problems$^{3-7}$ involve the use of the Frenet-Serret equations,
requiring numerical approximations or the use of helical arc segments (each
having constant $\tau$ and $\kappa$). \ Specifying $\kappa$ and $\theta$ to
generate a solution may be useful if $\tau$ is not initially known. \ The
torsion $\tau$ can then be determined from $\theta$.

\section{Mathematical Development}

A local coordinate system having the property $\boldsymbol{T}=\boldsymbol{i}%
^{\prime}$ (figure 1) supports the definition of $\boldsymbol{N}$:%

\begin{equation}
\boldsymbol{N}=\boldsymbol{j}^{\prime}\cos\theta+\boldsymbol{k}^{\prime}%
\sin\theta.
\end{equation}
%

\begin{figure}
[ptb]
\begin{center}
\includegraphics[
height=4.2618in,
width=3.2647in
]%
{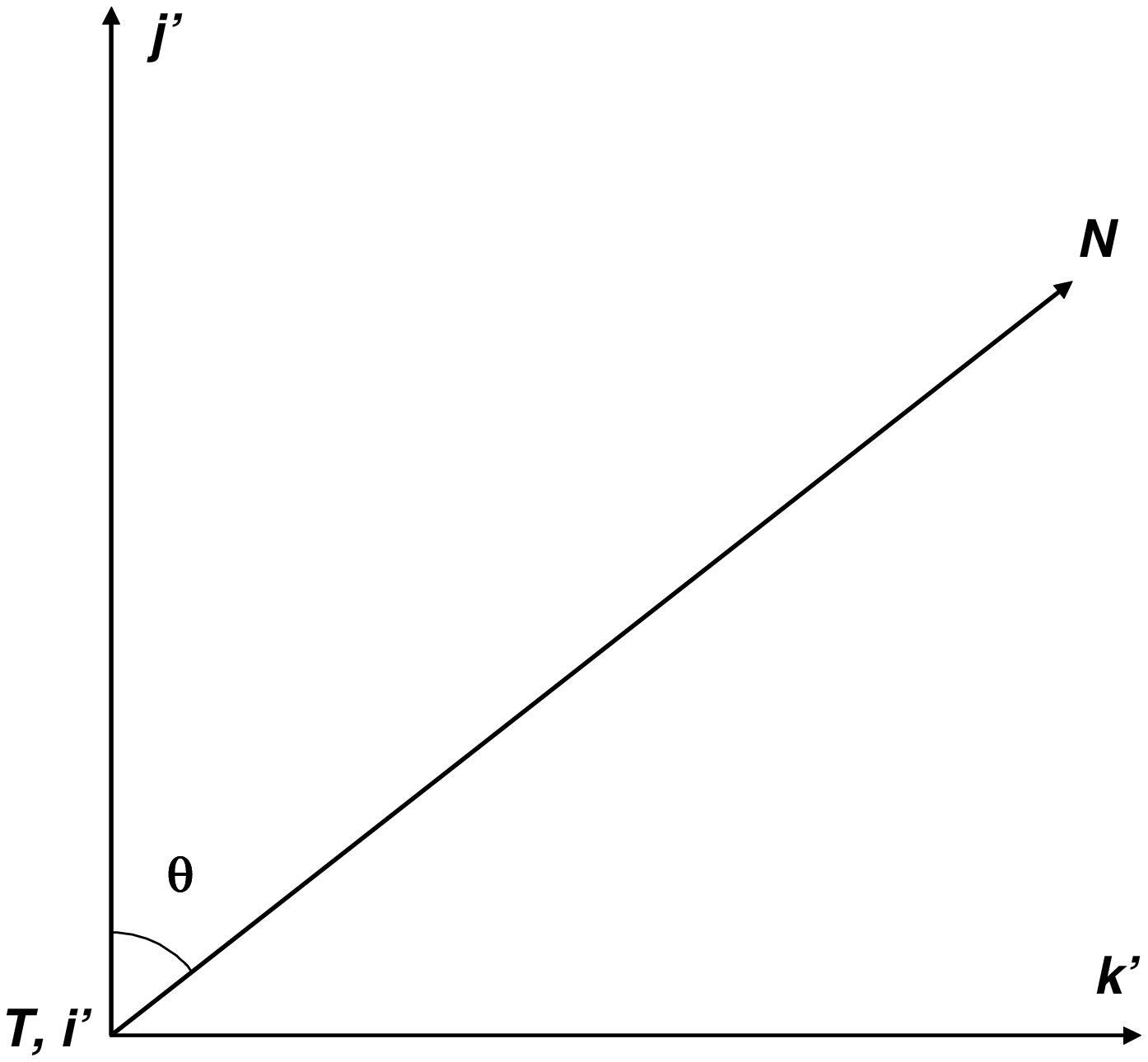}%
\caption{Local coordinate system; $\boldsymbol{T}$ is normal to the plane
containing $\boldsymbol{j}^{\prime}$ \& $\boldsymbol{k}^{\prime}$.}%
\end{center}
\end{figure}

The curvature $\kappa$ and the angle of rotation $\theta$ of the osculating
plane (containing $\boldsymbol{N}$ and $\boldsymbol{T}$) characterize the
curve. \ When the plane containing $\boldsymbol{T}$ and the global coordinate
$\boldsymbol{j}$ is normal to $\boldsymbol{k}^{\prime}$ (figure 2) then%

\begin{equation}
\boldsymbol{k}^{\prime}=\frac{\boldsymbol{T}\times\boldsymbol{j}}{\left\vert
\boldsymbol{T}\times\boldsymbol{j}\right\vert }=\frac{-\boldsymbol{i}%
T_{k}+\boldsymbol{k}T_{i}}{\sqrt{1-T_{j}^{2}}}.
\end{equation}
%

\begin{figure}
[ptb]
\begin{center}
\includegraphics[
trim=0.000000in 0.000000in 0.000000in 1.797346in,
height=3.2967in,
width=3.3451in
]%
{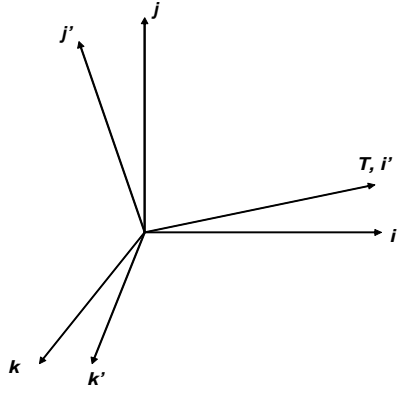}%
\caption{Angular orientation of the local coordinate system with respect to
the global coordinate system.}%
\end{center}
\end{figure}

Equation (6) breaks down when $\boldsymbol{T}=\pm\boldsymbol{j}$, requiring an
alternate expression for $\boldsymbol{k}^{\prime}$ (developed in section 4).
\ However, when $\boldsymbol{T}\neq\pm\boldsymbol{j}$,%

\begin{equation}
\boldsymbol{j}^{\prime}=\boldsymbol{k}^{\prime}\times\boldsymbol{T}%
=\frac{-\boldsymbol{i}T_{i}T_{j}+\boldsymbol{j}\left(  1-T_{j}^{2}\right)
-\boldsymbol{k}T_{j}T_{k}}{\sqrt{1-T_{j}^{2}}}\text{.}%
\end{equation}

Substituting (7) and (6) into (5):%

\begin{align}
N_{i}  &  =\frac{1}{\kappa}\frac{dT_{i}}{ds}=\frac{-T_{k}\sin\theta-T_{i}%
T_{j}\cos\theta}{\sqrt{1-T_{j}^{2}}};\\
N_{j}  &  =\frac{1}{\kappa}\frac{dT_{j}}{ds}=\cos\theta\sqrt{1-T_{j}^{2}};\\
N_{k}  &  =\frac{1}{\kappa}\frac{dT_{k}}{ds}=\frac{T_{i}\sin\theta-T_{j}%
T_{k}\cos\theta}{\sqrt{1-T_{j}^{2}}}.
\end{align}

Equation (9) can be integrated directly:%

\begin{equation}
\int_{T_{j0}}^{T_{j}}\frac{dT_{j}}{\sqrt{1-T_{j}^{2}}}=\int_{s_{0}}^{s}%
\kappa\cos\theta d\sigma;
\end{equation}

leading to%

\begin{align}
\sin^{-1}T_{j}  &  =\sin^{-1}T_{j0}+\int_{s_{0}}^{s}\kappa\cos\theta
d\sigma;\nonumber\\
T_{j}  &  =\sin\left[  \sin^{-1}T_{j0}+\int_{s_{0}}^{s}\kappa\cos\theta
d\sigma\right]  =\sin\delta;\nonumber\\
T_{j}  &  =T_{j0}\cos\int_{s_{0}}^{s}\kappa\cos\theta d\sigma+\sqrt
{1-T_{j0}^{2}}\sin\int_{s_{0}}^{s}\kappa\cos\theta d\sigma.
\end{align}

Equation (8) is solved by noting that $T_{k}=\sqrt{1-T_{j}^{2}-T_{i}^{2}%
}=\sqrt{\cos^{2}\delta-T_{i}^{2}}$ and introducing the variable $\beta$ so that%

\begin{align}
T_{i}  &  =\cos\delta\cos\beta;\\
T_{k}  &  =\cos\delta\sin\beta.
\end{align}

Substituting into (8):%

\begin{align}
\frac{dT_{i}}{ds}  &  =-\kappa\cos\theta\sin\delta\cos\beta-\cos\delta
\sin\beta\frac{d\beta}{ds}\nonumber\\
&  =\frac{-\kappa\sin\theta\cos\delta\sin\beta-\kappa\cos\theta\cos\delta
\cos\beta\sin\delta}{\cos\delta}\nonumber\\
&  =-\kappa\sin\theta\sin\beta-\kappa\cos\theta\cos\beta\sin\delta.
\end{align}

Equation (15) simplifies to%

\begin{equation}
\frac{d\beta}{ds}=\frac{\kappa\sin\theta}{\cos\delta};
\end{equation}

or%

\begin{equation}
\beta=\beta_{0}+\int_{s_{0}}^{s}\frac{\kappa\sin\theta}{\cos\delta}%
d\sigma=\cos^{-1}\left(  \frac{T_{i}}{\cos\delta}\right)  ;
\end{equation}

\bigskip so that%

\begin{equation}
T_{i}=T_{i_{0}}\frac{\cos\delta}{\cos\delta_{0}}\cos\int_{s_{0}}^{s}%
\frac{\kappa\sin\theta}{\cos\delta}d\sigma-T_{k_{0}}\frac{\cos\delta}%
{\cos\delta_{0}}\sin\int_{s_{0}}^{s}\frac{\kappa\sin\theta}{\cos\delta}%
d\sigma;
\end{equation}

where%

\begin{equation}
\cos\delta=\sqrt{1-T_{j0}^{2}}\cos\int_{s_{0}}^{s}\kappa\cos\theta
d\sigma-T_{j0}\sin\int_{s_{0}}^{s}\kappa\cos\theta d\sigma.
\end{equation}

The solution for $T_{k}$ follows from (14) and (16):%

\begin{equation}
T_{k}=T_{k_{0}}\frac{\cos\delta}{\cos\delta_{0}}\cos\int_{s_{0}}^{s}%
\frac{\kappa\sin\theta}{\cos\delta}d\sigma+T_{i_{0}}\frac{\cos\delta}%
{\cos\delta_{0}}\sin\int_{s_{0}}^{s}\frac{\kappa\sin\theta}{\cos\delta}%
d\sigma.
\end{equation}

It can be easily verified that (12), (18), and (20) meet the requirement:%

\begin{equation}
\kappa=\left\vert \frac{d\boldsymbol{T}}{ds}\right\vert .
\end{equation}

Generating an expression for the torsion $\tau$ requires first computing
$\boldsymbol{N}$ by substituting (12)-(14) into (8)-(10):%
\begin{align}
N_{i}  &  =-\cos\theta\sin\delta\cos\beta-\sin\beta\sin\theta;\\
N_{j}  &  =\cos\theta\cos\delta;\\
N_{k}  &  =-\cos\theta\sin\delta\sin\beta+\cos\beta\sin\theta.
\end{align}

Next, $\boldsymbol{B=T\times N}$:

\begin{align}
B_{i}  &  =\sin\delta\sin\theta\cos\beta-\cos\theta\sin\beta;\\
B_{j}  &  =-\cos\delta\sin\theta;\\
B_{k}  &  =\sin\delta\sin\theta\sin\beta+\cos\theta\cos\beta.
\end{align}

Equation (28) expresses the torsion as a function of $\theta$:%

\begin{equation}
\tau=\left\vert \frac{d\boldsymbol{B}}{ds}\right\vert =\frac{d\theta}%
{ds}-\kappa\tan\delta\sin\theta.
\end{equation}

Equation (29) expresses $\tau$ in terms of components of $\boldsymbol{T}$ and
$\boldsymbol{B}$:%

\begin{equation}
\tau=\frac{d\theta}{ds}+\frac{\kappa T_{j}B_{j}}{1-T_{j}^{2}}.
\end{equation}

Finally, (2) serves as a check on the solutions for $\boldsymbol{T}$,
$\boldsymbol{N}$, $\boldsymbol{B}$, and $\tau$.

\section{Discussion}

\bigskip Integrating (29) leads to the following expression for $\theta$:%

\begin{equation}
\theta=\theta_{0}+\int_{s_{0}}^{s}\left(  \tau-\frac{\kappa T_{j}B_{j}%
}{1-T_{j}^{2}}\right)  d\sigma.
\end{equation}

Equation (30) indicates a unique value of $\theta$ for each specified value of
$\tau$ when $T_{j}\neq\pm1$. \ Thus, (12), (18), and (20) indirectly solve (1)-(3).

The angle $\theta$ can also be expressed in terms of components of
$\boldsymbol{T}$, $\boldsymbol{N}$, $\boldsymbol{B}$:%

\begin{equation}
\theta=-\sin^{-1}\frac{B_{j}}{\sqrt{1-T_{j}^{2}}}=\cos^{-1}\frac{N_{j}}%
{\sqrt{1-T_{j}^{2}}}=-\tan^{-1}\frac{B_{j}}{N_{j}}.
\end{equation}

\subsection{Constant $\theta$}

An explicit solution often results when $\theta$ is constant. \ Setting
$T_{i0}=1$, so that $\beta_{0}=\delta_{0}=0$ (and setting $s_{0}=0$) leads to%

\begin{align}
\delta &  =\int_{0}^{s}\kappa(\sigma)\cos\theta_{0}d\sigma;\\
\beta &  =2\tan\theta_{0}\tanh^{-1}\left(  \tan\delta/2\right)  ;
\end{align}

so that%

\begin{align}
T_{i}  &  =\cos\left[  2\tan\theta_{0}\tanh^{-1}\left(  \tan\delta/2\right)
\right]  \cos\int_{0}^{s}\kappa(\sigma)\cos\theta_{0}d\sigma;\\
T_{j}  &  =\sin\int_{0}^{s}\kappa(\sigma)\cos\theta_{0}d\sigma;\\
T_{k}  &  =\sin\left[  2\tan\theta_{0}\tanh^{-1}\left(  \tan\delta/2\right)
\right]  \cos\int_{0}^{s}\kappa(\sigma)\cos\theta_{0}d\sigma.
\end{align}

The torsion becomes%

\begin{equation}
\tau(s)=-\kappa(s)\sin\theta_{0}\tan\int_{0}^{s}\kappa(\sigma)\cos\theta
_{0}d\sigma\text{.}%
\end{equation}

As an example, when%

\begin{equation}
\kappa=\kappa_{0}e^{-s^{2}}\text{,}%
\end{equation}

\begin{align}
T_{i}  &  =\cos\left[  \frac{\kappa_{0}\sqrt{\pi}}{2}\operatorname{erf}%
(s)\cos\theta_{0}\right]  \cos\left[  2\tan\theta_{0}\tanh^{-1}\left(
\tan\left[  \frac{\kappa_{0}\sqrt{\pi}}{4}\operatorname{erf}(s)\cos\theta
_{0}\right]  \right)  \right]  \text{;}\\
T_{j}  &  =\sin\left[  \frac{\kappa_{0}\sqrt{\pi}}{2}\operatorname{erf}%
(s)\cos\theta_{0}\right]  \text{;}\\
T_{k}  &  =\cos\left[  \frac{\kappa_{0}\sqrt{\pi}}{2}\operatorname{erf}%
(s)\cos\theta_{0}\right]  \sin\left[  2\tan\theta_{0}\tanh^{-1}\left(
\tan\left[  \frac{\kappa_{0}\sqrt{\pi}}{4}\operatorname{erf}(s)\cos\theta
_{0}\right]  \right)  \right]  \text{;}%
\end{align}

\begin{equation}
\tau(s)=-\kappa_{0}e^{-s^{2}}\sin\theta_{0}\tan\left[  \frac{\kappa_{0}%
\sqrt{\pi}}{2}\operatorname{erf}(s)\cos\theta_{0}\right]  \text{.}%
\end{equation}

\subsection{Constant $\kappa$}

When $\kappa=\kappa_{0}$ but $\theta\neq\theta_{0}$, the solution will
typically involve undetermined integrals. \ For example, when $\kappa
=\kappa_{0}$ and $\theta=\kappa_{0}s$,%

\begin{align}
T_{i}  &  =\cos\left(  \sin\kappa_{0}s\right)  \cos\int_{0}^{s}\frac
{\kappa_{0}\sin\kappa_{0}\sigma}{\cos\left(  \sin\kappa_{0}\sigma\right)
}d\sigma;\\
T_{j}  &  =\sin\left(  \sin\kappa_{0}s\right)  ;\\
T_{k}  &  =\cos\left(  \sin\kappa_{0}s\right)  \sin\int_{0}^{s}\frac
{\kappa_{0}\sin\kappa_{0}\sigma}{\cos\left(  \sin\kappa_{0}\sigma\right)
}d\sigma;
\end{align}

and%

\begin{equation}
\tau(s)=\kappa_{0}-\kappa_{0}\tan\left(  \sin\kappa_{0}s\right)  \sin
\kappa_{0}s.
\end{equation}

\subsection{\bigskip Constant $\kappa$ and $\theta$}

When $\kappa=\kappa_{0}$ and $\theta=\theta_{0}$, (34)-(37) become:%

\begin{align}
T_{i}  &  =\cos\left(  \kappa_{0}s\cos\theta_{0}\right)  \cos\left[
2\tan\theta_{0}\tanh^{-1}\left(  \tan\left[  \frac{\kappa_{0}s}{2}\cos
\theta_{0}\right]  \right)  \right]  ;\\
T_{j}  &  =\sin\left(  \kappa_{0}s\cos\theta_{0}\right)  ;\\
T_{k}  &  =\cos\left(  \kappa_{0}s\cos\theta_{0}\right)  \sin\left[
2\tan\theta_{0}\tanh^{-1}\left(  \tan\left[  \frac{\kappa_{0}s}{2}\cos
\theta_{0}\right]  \right)  \right]  ;
\end{align}

\begin{equation}
\tau(s)=-\kappa_{0}\sin\theta_{0}\tan\left(  \kappa_{0}s\cos\theta_{0}\right)
.
\end{equation}

When $\theta_{0}=\pi/2$, $\tau(s)=0$, confining $\boldsymbol{T}$ and
$\boldsymbol{N}$ to a plane. \ When $\boldsymbol{T}$ aligns with
$\boldsymbol{j}$, $\tau\rightarrow\infty$ in (50), and the equations break down.

\section{Alternate Set of Equations}

The equations break down when $T_{j}\rightarrow\pm1$, requiring a different
orientation for the local coordinate system. \ The angle of rotation of the
osculating plane is designated $\phi$ here. In general, $\phi\neq\theta$,
reflecting differences in angular orientation between the local and global
coordinate systems for the two cases. \ Defining $\boldsymbol{k}^{\prime}$ as
the normal to the plane containing $\boldsymbol{T}$ and $\boldsymbol{i}$, i.e.,%

\begin{equation}
\boldsymbol{k}^{\prime}=\frac{\boldsymbol{i}\times\boldsymbol{T}}{\left\vert
\boldsymbol{i}\times\boldsymbol{T}\right\vert }=\frac{-\boldsymbol{j}%
T_{k}+\boldsymbol{k}T_{j}}{\sqrt{1-T_{i}^{2}}}.
\end{equation}

The $\boldsymbol{j}^{\prime}$ unit vector becomes:%

\begin{equation}
\boldsymbol{j}^{\prime}=\boldsymbol{k}^{\prime}\times\boldsymbol{T}%
=\frac{-\boldsymbol{i}\left(  1-T_{i}^{2}\right)  +\boldsymbol{j}T_{i}%
T_{j}+\boldsymbol{k}T_{i}T_{k}}{\sqrt{1-T_{i}^{2}}}.
\end{equation}

Substituting into the expression for $\boldsymbol{N}$:%

\begin{align}
N_{i}  &  =\frac{1}{\kappa}\frac{dT_{i}}{ds}=-\cos\phi\sqrt{1-T_{i}^{2}};\\
N_{j}  &  =\frac{1}{\kappa}\frac{dT_{j}}{ds}=\frac{-T_{k}\sin\phi+T_{i}%
T_{j}\cos\phi}{\sqrt{1-T_{i}^{2}}};\\
N_{k}  &  =\frac{1}{\kappa}\frac{dT_{k}}{ds}=\frac{T_{j}\sin\phi+T_{i}%
T_{k}\cos\phi}{\sqrt{1-T_{i}^{2}}}.
\end{align}

Equations (53)-(55) have the following solution:%

\begin{align}
T_{i}  &  =\sin\gamma;\\
T_{j}  &  =\cos\gamma\cos\alpha;\\
T_{k}  &  =\cos\gamma\sin\alpha;\\
N_{i}  &  =-\cos\gamma\cos\phi;\\
N_{j}  &  =\sin\gamma\cos\alpha\cos\phi-\sin\alpha\sin\phi;\\
N_{k}  &  =\sin\gamma\sin\alpha\cos\phi+\cos\alpha\sin\phi;\\
B_{i}  &  =\cos\gamma\sin\phi;\\
B_{j}  &  =\sin\gamma\cos\alpha\sin\phi+\sin\alpha\cos\phi;\\
B_{k}  &  =-\sin\gamma\sin\alpha\sin\phi+\cos\alpha\cos\phi;\\
\tau &  =\frac{d\phi}{ds}-\kappa\tan\gamma\sin\phi=\frac{d\phi}{ds}%
-\frac{\kappa B_{i}T_{i}}{1-T_{i}^{2}}.
\end{align}

Here%

\begin{align}
\gamma &  =\sin^{-1}T_{i0}-\int_{s_{0}}^{s}\kappa\cos\phi d\sigma;\\
\alpha &  =\cos^{-1}\left(  \frac{T_{j}}{\cos\gamma}\right)  =\alpha_{0}%
+\int_{s_{0}}^{s}\frac{\kappa\sin\phi}{\cos\gamma}d\sigma;\\
T_{j}  &  =T_{j_{0}}\frac{\cos\gamma}{\cos\gamma_{0}}\cos\int_{s_{0}}^{s}%
\frac{\kappa\sin\phi}{\cos\gamma}d\sigma-T_{k_{0}}\frac{\cos\gamma}{\cos
\gamma_{0}}\sin\int_{s_{0}}^{s}\frac{\kappa\sin\phi}{\cos\gamma}d\sigma;\\
\phi &  =-\tan^{-1}\frac{B_{i}}{N_{i}}.
\end{align}

Even though $\theta$ and $\phi$ both represent the angle of rotation of the
osculating plane, (31) and (69) differ because of differences in angular
orientation of the local coordinate system.

When $T_{j}\rightarrow\pm1$ or $T_{i}\rightarrow\pm1$, switching from one set
of equations to another avoids numerical difficulties.

\section{Concluding Remarks}

Unlike the Frenet-Serret equations, (8)-(10) are nonlinear, and do not involve
$\boldsymbol{N}$, $\boldsymbol{B}$, or $\tau$. \ The solution (in terms of
$\kappa$ and $\theta$) indirectly solves the Frenet-Serret equations, and
leads to a precise definition of $\tau$ as a function of $\kappa$ and $\theta
$. \ A unique value of $\theta$ can be obtained for each specified value of
$\tau$ through a first order ordinary differential equation. \ The equations
break down when $\boldsymbol{T}\rightarrow\pm\boldsymbol{j}$, requiring an
alternative set of equations that break down when $\boldsymbol{T}%
\rightarrow\pm\boldsymbol{i}$. \ The expressions for the angle of the
osculating plane in the two approaches differ because of differences in the
angular orientation of the local coordinate system.

\textbf{Acknowledgement} \ \textit{This work was funded by the Office of Naval
Research, Code 321US (M. Vaccaro).}

\section{References}

\begin{enumerate}
\item M. P. do Carmo (1976). \textit{Differential Geometry of Curves and
Surfaces}. \ Prentice-Hall, Englewood Cliff, NJ.

\item B. Divjak (1997). \textit{Mathematical Communications} \textbf{2}, 143-147.

\item K. Nakayama, H. Segur, \& M. Wadati (1992). \ \textit{Phys. Rev. Lett.}
\textbf{69}, 2603-2606.

\item Y. Kats, D. A. Kessler, \& Y. Rabin (2002). \ \textit{Phys. Rev. E}
\textbf{65}, 020801(R).

\item H. Hasimoto (1972). \ \textit{J. Fluid Mech.} \textbf{51}, 477-485.

\item G. Arreaga-Garcia, H. Villegas-Brena, \& J. Saucedo-Morales (2004).
\textit{J. Phys. A: \ Math. Gen.} \textbf{37}, 9419-9438.

\item A. C. Hausrath \& A. Goriely (2006). \ \textit{Protein Science}
\textbf{15}, 753-760.
\end{enumerate}

\end{document}